\documentclass[a4paper,11pt]{article}

\newtheorem{theorem}{Theorem}
\newtheorem{proposition}[theorem]{Proposition}
\newtheorem{lemma}[theorem]{Lemma}
\newtheorem{counter-example}[theorem]{Counter-example}

\newtheorem{definition}[theorem]{Definition}

\newtheorem{remark}[theorem]{Remark}

\newcommand{\qed}{\nobreak \ifvmode \relax \else
      \ifdim\lastskip<1.5em \hskip-\lastskip
      \hskip1.5em plus0em minus0.5em \fi \nobreak
      \vrule height0.75em width0.5em depth0.25em\fi}

\makeindex 

\usepackage{makeidx}     
\usepackage{cite}        
\usepackage{amsfonts}
\usepackage{amssymb}
\usepackage{amscd}
\usepackage{epsfig}
\usepackage{amsrefs}
\usepackage{url}
\usepackage[all]{xy}
\usepackage{pb-diagram}

\def\C{\mathbf{C}}
\def\F{\mathbf{F}}

\def\rmi{\uppercase\expandafter{\romannumeral1}}
\def\rmii{\uppercase\expandafter{\romannumeral2}}
\def\rmiii{\uppercase\expandafter{\romannumeral3}}
\def\rmv{\uppercase\expandafter{\romannumeral5}}
\def\rmvi{\uppercase\expandafter{\romannumeral6}}
\def\rmvii{\uppercase\expandafter{\romannumeral7}}
\def\rmviii{\uppercase\expandafter{\romannumeral8}}

\def\U{\mathcal{U}}

\def\can#1{\left\langle#1\right\rangle}

\def\PB{\left\{\cdot\,,\cdot\right\}}

\def\pb#1{\left\{#1\right\}}

\def\lb#1{\[#1\]}

\def\({\left(}
\def\){\right)}
\def\[{\left[}
\def\]{\right]}

\def\diff{\mathsf{d}}

\def\Ker{\mathop{\rm Ker}\nolimits}

\catcode`\é=13 \defé{\'e}
\catcode`\è=13 \defè{\`e}
\catcode`\ê=13 \defê{\^e}

\newenvironment{equation*}[1][1.5]
  {$$\renewcommand{\arraystretch}{#1}
      \begin{array}{rcl}}
      {\end{array}$$}

\def\comment#1{}  

\def\p{\partial}
\def\pp#1#2{\frac{\p #1}{\p #2}}

\begin{document} 
\nocite{*}

\title{Cotangent paths as coisotropic subsets for local functions}

\author{
  Camille Laurent-Gengoux\\
  \texttt{camille.laurent-gengoux@univ-lorraine.fr}
\and
  Yahya Turki\\
  \texttt{yahya.turki@univ-lorraine.fr}
  \thanks{Universit\'e de Lorraine, 
Institut Elie Cartan de Lorraine
UMR 7502, Metz, F-57045, France.
}
\thanks{Universit\'e de Monastir, Facult\'e des Sciences de Monastir, Avenue de l'Environnement 5019 Monastir, Tunisie.
}
}
                    
\maketitle
\begin{abstract}
We establish a local function version of a classical result claiming that a bivector field on a manifold $M$ is Poisson if and only if cotangent paths form a coisotropic set
of the infinite dimensional symplectic manifold of paths valued in $T^*M$.
Our purpose here is to prove this result without using the Banach manifold setting, setting which fails in the periodic case because cotangent loops do not 
form a Banach sub-manifold. Instead, we use local functions on the path space, a point of view that allows to speak of a coisotropic set.
\end{abstract}

\tableofcontents
\newpage

\section{Introduction}

For every symplectic manifold $S$, both loop and path spaces with values in $S$ comes equipped with natural symplectic manifold structures (certainly, infinite dimensional, 
which raises concern about the meaning of "symplectic" in that context, but we can make sense of this affirmation by using local functions).
Since the cotangent space $T^*M$ of a manifold $M$ comes equipped with a symplectic form, 
paths and loops valued in $T^{*}M$ are in particular a symplectic manifold. 
In this section, our purpose is give a definition making sense, with local functions:
\begin{enumerate}
\item of a Lie bracket on local functions on loops or paths valued in $T^*M$, 
\item of what we mean by coisotropic subset of the space of loops or paths in $T^{*}M$. 
\end{enumerate}
Then we intend to show that a bivector field is Poisson if and only if a certain subset (called cotangent paths) of the space of paths in $T^*M$ is coisotropic.

Let us say a few words on the history of this result, which is now a classical result of Poisson Sigma models. 

To our knowledge, it was first establish by Shaller and Strobl in Equation (10) of \cite{SchallerStrobl}, generalized and implicitly 
stated as an equivalence by Klim{\v{c}}{\'{\i}}k and Strobl in Equation (7) in \cite{strobelklim}. More precisely,
the authors of \cite{strobelklim} prove  a more general result for WZW manifolds, also called Poisson
manifolds with background or twisted Poisson manifolds \cite{KosmannSchwarzbachYvette}; also, the authors of \cite{strobelklim} do not use the word coisotropic,
but simply give functions whose zero locus are cotangent paths and prove that their 
Poisson bracket is zero on this zero locus if and only if the (twisted) Jacobi identity holds.  Of course, this computation is also 
the crucial point of the proof presented in this article, and, indeed, we shall 
compute and exploit the same very brackets later on. Note that, strictly speaking, this computation is not sufficient to show that cotangent paths form a coisotropic sub-manifold.  
The ideal generated by these functions is a priori smaller than the ideal of all functions vanishing on cotangent space, so we can not, using Equation (7) 
in \cite{strobelklim}, exclude that there might be local functions vanishing on all
cotangent paths whose bracket does not vanish when evaluated at a cotangent path.

Later on, this result was made mathematically very precise and exploited fruitfully by Cattaneo and Felder who used $C^1$-paths valued in $T^*M$, a set that comes equipped with a Banach manifold structure. 
In  \cite{CattaneoFelderCoisotropic}, the result, stated with the language of coisotropic subset, was announced.
A complete proof appeared in \cite{Cattaneo}. Using the Banach manifold setting, they show that  coisotropic paths do
form a sub-Banach manifold of the Banach manifold of $C^1$-paths valued in $T^*M$ which happens to be coisotropic if and only if $\pi$ is Poisson. 
Unfortunately the arguments presented in \cite{Cattaneo} fails for loops (see Section \ref{sec:counterExample} below). 
The issue may seem quite technical: cotangent loops do not form a submanifold, although is it certainly coisotropic, for $\pi$ Poisson, at each cotangent loop around which 
cotangent loops do form a submanifold. 

This result, at least the "if" part, is used in \cite{CattaneoFelderForK} to reconstruct a star product on the initial Poisson manifold.
It is also shown in \cite{CattaneoFelder} that
the symplectic groupoid, if it exists, is obtained by symplectic reduction of this coisotropic submanifold. 
This theorem is also extended to Courant algebras in \cite{AlekseevStrobl}.
Also, in \cite{Yahya}, one the author of this article has intended to characterize Poisson bivectors among bivector fields through the properties of some Lagrangians:
for instance, \cite{Yahya} introduces a Lagrangian, inspired by Poisson $\sigma$-models, whose stationary points are quasi-Hamiltonians 
(which is slightly weaker than cotangent paths) if and only if $\pi$ is Poisson.

Our purpose here is to prove or reprove the characterization of Poisson bivectors among bivector fields
whose history has been described above by working with local functions on the set of smooth paths, for which the Banach manifold setting can be avoided. 
By doing so, we claim that we are able to deal with periodic paths and non-periodic paths altogether, while the Banach manifold setting only 
works for non-periodic paths. 

In this paper, we only consider the case where $M$ is an open subset of ${\mathbb R}^n$, which does not lead to any loss of generality.
The paper is organized as follows. Section \ref{sec:local} is devoted to the study of local functionals on two subsetsof the set of all smooth paths, namely semi-free paths and loops, and to the 
definition of a Lie bracket on those, see Equation (\ref{eq:Liebracket}). In the course of Section \ref{sec:local}, coisotropic subsets are defined
(and we insist that we define coisotropic subsets, not coisotropic submanifolds).
The main objective of this paper is to show Theorem \ref{theo:1} (resp. \ref{theo:2}) which states that a  bivector field is Poisson if and only if semi-free (resp. periodic) 
cotangent paths form a coisotropic subset of the set of semi-free (resp.  periodic) paths in $T^{*}M$.
In Section \ref{Non-periodic case}, we prove Theorem \ref{theo:1} (for 
 semi free paths).
Section \ref{periodiccase} deals with the periodic case and we prove Theorem \ref{theo:2} in the first Section \ref{lecasperiodique}. In Section \ref{sec:counterExample},
we give an example showing that $C_{\pi}^{S^1}(M)$ is not a Banach manifold in general.

\section{Local functionals on subsets of paths}

\label{sec:local}

The space ${\mathfrak P}$ of smooth paths valued in a manifold $N$ is a Fréchet manifold \cite{PM}, and indeed inherits a topology modeled on the Fr\'echet topology,
but we intend to avoid this infinite manifold setting by dealing with local functions (also called local functionals).
Those functions can be defined for arbitrary source and target manifolds (see Equations (8-9)
\cite{Khavkine} for instance, or Theorem 10.4 in \cite{PM}), but we are only interested in the case $I= \lb{0,1}$:

\begin{definition}\label{def:fonctions_locales_general}
A ${\mathbb R}$-valued function $F$ on ${\mathfrak P}$ is called local when there exists,
an integer $k \in {\mathbb N}$ and a ${\mathbb R}$-valued  smooth function $f$ from the manifold $J^kN$ of $k$-jets of maps from ${\mathbb R}$ to $N$:
$$ \begin{array}{rrcl}
F_f :& {\mathfrak P} & \to & {\mathbb R}\\  
 & a &\longmapsto &  \int_{I} f \left( J^k a(t) \right) {\rm d}t,  
\end{array}
$$
with $J^k a(t) \in J^k N$ being the $k$-th jet of $a$ at $t$. We denote by $\F({\mathfrak P})$ the vector space of local functionals.
\end{definition}

When $N$ is an open subset of ${\mathbb R}^n$, the definition can be simplified:

\begin{definition}\label{def:fonctions_locales_N}
A ${\mathbb R}$-valued function $F$ on ${\mathfrak P}$ is called local when there exists, an integer
$k \in {\mathbb N}$ and a ${\mathbb R}$-valued  smooth function $f$ from $N \times {\mathbb R}^{kn}$
such that
$F=F_f$ where
$$ \begin{array}{rrcl}
F_f :& P & \to & {\mathbb R}\\  
 & a &\longmapsto &  \int_{I} f \left(t,a(t),a'(t),\dots,a^{(k)}(t) \right) {\rm d}t.  
 \end{array}
$$
\end{definition}

\begin{remark}\label{rmk:non-unicite}
For a local functional $F$ the function $f$ such that $F_f = F$ is not unique.
For example, when $N={\mathbb R}$, consider a function  $f$ of the form:
\begin{eqnarray*} f &= & \frac{\diff g \left(q,q',\dots,q^{(k)},t \right)  }{\diff t} \\
&=&  \frac{\partial g \left(q,q',\dots,q^{(k)},t \right)  }{\partial t} + \sum_{i \geq 0} q^{(i+1)}\frac{\partial g \left(q,q',\dots,q^{(k)},t \right)  }{\partial q^{(i)}}, \end{eqnarray*}
with $g\left(q,q',\dots,q^{(k)},p,p',\dots,p^{(l)},t \right)=0$ for $t=0$ and $t=1$.
It is clear that $F_f= 0$ although $f \neq 0$.
\end{remark}

One can enlarge slightly Definition \ref{def:fonctions_locales_general} and call local functions
functions on  ${\mathfrak P}$ such that every element of  ${\mathfrak P}$ admits a neighborhood in which they are of the form $F_f$ for some function $f$ defined on the space of jets,
see \cite{BDLR}. We do not need this extension at the moment. We now define local functions on subsets (not necessarily submanifolds in any sense) of ${\mathfrak P}$.

\begin{definition}
A function on a subset $\tilde{\mathfrak P}$ of ${\mathfrak P}$ is called \emph{local} when it is the restriction to this subset of a local function on ${\mathfrak P}$. We denote by $\F(\tilde{\mathfrak P})$ the vector space of local functions on $\tilde{\mathfrak P}$. 
\end{definition}

We are mainly interested on two subsets of the set of smooth paths: periodic smooth paths and semi-free smooth paths.
\emph{Periodic smooth paths} or \emph{loops} are just those paths whose values and $k$-th order derivatives at $t=0$ and $t=1$ coincide for all $k\geq 0$. Loops can be seen 
as smooth maps from $S^1$ to $N$ and shall be denoted by ${\mathfrak L}$. When $N=E \to M$ is a vector bundle, we call \emph{semi-free paths} and denote by $\tilde{\mathfrak P}$ 
the smooth paths $a \in {\mathfrak P}$ valued in $E$ that
\begin{enumerate}
\item[1)] starts from the zero section i.e. $a(0)=0 \in E_{q(0)}$ (with $q:I \to M$ being the base path of $a$, i.e. the image of $a$ through the natural projection $E \to M$),
\item[2)] arrives at zero section i.e. $a(1)=0 \in E_{q(1)}$,
\item[3)] whose $k$-th order derivatives at $t=0$ and $t=1$ are all equal to zero for all $k \geq 1$.
\end{enumerate}
\begin{center}
\texttt{ From now, and until the end of this section,  $N$ is $T^{*}M $ with $M$ an open subset of ${\mathbb R}^{n}$.}
\end{center}
We identify $T^{*}M$ with $M  \times  {\mathbb R}^{n}$, and the canonical projection $T^*M \to M $ (of the cotangent bundle onto its base) with the projection onto the first component: 
 $$(q,p) \to q  \quad \quad \hbox{ for all } q,p \in{\mathbb R}^{n}.$$
 The isomorphism $T^{*}M \simeq M  \times  {\mathbb R}^{n}$ allows to write $a=(q,p)$ where $q$ and $p$ are paths in $ M$ and in $ {\mathbb R}^{n}$ respectively. The path $t \mapsto q(t)$, which is the projection of $a$ on the base manifold $M$, is the \emph{base path}.

A function $F$ on ${\mathfrak P}$ with values in ${\mathbb R}$ is therefore local when there exist integers
$k,l \in {\mathbb N}$ and a smooth function $f$ from $M \times I \times {\mathbb R}^{(k+l)n}$ to ${\mathbb R} $
such that $F=F_f$ where
$$ \begin{array}{rrcl}
F_f :&\mathfrak P & \to & {\mathbb R}\\  
 & a =(q,p)&\longmapsto &  \int_{I} f \left(t,q(t),q'(t),\dots,q^{(k)}(t),p(t),p'(t),\dots,p^{(l)}(t) \right) {\rm d}t.  
\end{array}
$$
Also, in view of the definition above, semi-free paths $\tilde{{\mathfrak P}} $ are those paths $a = (q,p)\in {\mathfrak P}$ such that $p(0)=0$,  $p(1)=0$,
and the $k$-th order derivatives of $q$ and $p$ vanish at $t=0$ and $t=1$ for all $k \geq 1$.

We now intend to define the following two items:
\begin{enumerate}
\item[1)] the gradient of a local function,
\item[2)] the Poisson bracket of two local functions. 
It is important to notice that local functions  $ \F(\tilde{\mathfrak P})$ on semi-free paths 
(or local functions $ \F({\mathfrak L})$ on loops) do not form an algebra, since it is not stable by product, 
hence this Poisson bracket is in fact a Lie bracket: we keep this name, however, because it extends to a Poisson bracket in the 
algebra generated by $\F(\tilde{\mathfrak P})$ or $\F({\mathfrak L})$.
\end{enumerate}

 We begin by defining the gradient of a local functional $F \in \F(\tilde{{\mathfrak P}})$ on ${\mathfrak P}$ at a semi-free path  or a loop $a =(q,p)\in \tilde{\mathfrak P}$.
By \emph{a vector tangent to ${\mathfrak P}$ at a point $a\in {\mathfrak P}$}, we simply mean a function $\delta_a  \in C^\infty (I , {\mathbb R}^n \times {\mathbb R^n})$. In view of the decomposition $T^*M = M \times {\mathbb R}^n$, we write
$ \delta_a =(\delta_q,\delta_p)$ where  $\delta_q,\delta_p \in C^\infty (I , {\mathbb R}^n)$. \emph{A vector tangent to $\tilde{\mathfrak P}$ at $ a \in  \tilde{\mathfrak P} $} is a vector $\delta_a=(\delta_q,\delta_p)$ tangent to $P$ such that $\delta_q'(t), \delta_p(t)$ vanish together with all their derivatives at $t=0$ and $t=1$. \emph{A vector tangent to ${\mathfrak L}$ at $a \in {\mathfrak L} $} is a vector tangent to $a \in {\mathfrak P}$ such that $\delta_a \in  {\mathcal C}^\infty (S^1 , {\mathbb R}^n \times {\mathbb R^n})$.

This name of tangent vector comes from the easily checked property that for all $\epsilon \in {\mathbb R}$ small enough, the path $t \mapsto a(t)+\epsilon \delta_a(t)$ is an element of ${\mathfrak P},\tilde{\mathfrak P}, {\mathfrak L}$ (depending on the context) and $\delta_a$ is the derivative at $\epsilon=0$ of this expression\footnote{That is, the tangent vector of ${\mathfrak P}$, $\tilde{\mathfrak P}$ or $ {\mathfrak L}$ coincides with the tangent cone.}.

\begin{lemma}  \label{gradient-semi}
For any local function $F$ on $\tilde{\mathfrak P}$, there exists one and only one element 
$\nabla _a F = (A_F,B_F,\alpha^0_F ,\alpha^1_F ) \in \C^{\infty}(I,{\mathbb R}^{n}\times  {\mathbb R}^{n}) \times  ({\mathbb R}^n \times {\mathbb R^n})$ satisfying that  for any tangent vector $\delta_a =(\delta_q,\delta_p)$ of $\tilde{\mathfrak P}$:
\begin{eqnarray*}
\left. \frac{\diff }{\diff \epsilon}\right|_{\epsilon=0} F(a+\epsilon \delta_a )
  &=   &\int_{I} \left\langle A_F ,\delta_q  \right\rangle {\rm d}t + \int_{I}  \left\langle B_F ,\delta_p  \right\rangle {\rm d}t \\ 
   &+   & \left\langle \alpha^1_F ,\delta_{q}(1)  \right\rangle - \left\langle \alpha^0_F ,\delta_{q}(0)  \right\rangle,
   \end{eqnarray*}   
where $\langle . , . \rangle$ is the canonical scalar product of ${\mathbb R}^n$.
We call this quadruple the gradient of $F$ at $a$. 
 \end{lemma}
 \begin{proof}
 The gradient is clearly unique if it exists because the restriction to $]-\epsilon ,1-\epsilon[$ of $ \delta_q$ and $\delta_p$ can  be chosen to be arbitrary smooth functions.
 
For the existence and an explicit formula for $ \nabla _a F $ choose a function  $f$  as in the definition (\ref{def:fonctions_locales_N}) such that $F=F_f$.
 By definition of a local function, we have for all small enough $\epsilon$:
$$F(a+\epsilon \delta _a) = \int_{I} f \left(q+\epsilon \delta _q,\dots,q^{(k)}+\epsilon \delta^{(k)} _q,p+\epsilon \delta _p,\dots,p^{(l)}+\epsilon \delta^{(l)} _p,t \right) {\rm d}t.$$
We differentiate with respect to $\epsilon$ and obtain, by denoting by $\frac{\partial f}{\partial q^{(i)}}$ or $ \frac{\partial f}{\partial p^{(j)}} $ 
the gradients of the function $f$ with respect the variable $ q^{(i)}$ or $p^{(j)}$
 \begin{eqnarray*}
& & \left. \frac{\diff F(a+\epsilon \delta _a)}{\diff \epsilon} \right|_{\epsilon=0} \\
  &=  &\int_{I} \left( \left< \frac{\partial f}{\partial q},\delta _q \right> + \left<\frac{\partial f}{\partial q'},\delta' _q \right> +\dots +\left< \frac{\partial f}{\partial q^{(k)}},
  \delta^{(k)} _q(t) \right> \right. \\  
  &+   & \left. \left<\frac{\partial f}{\partial p},\delta _p \right> +\left<\frac{\partial f}{\partial p'},\delta' _p  \right> + 
  \dots +\left<\frac{\partial f}{\partial p^{(l)}},\delta^{(l)} _p(t)  \right> \right)  {\rm d}t \\
  &=   & \int_{I} \left< \left. \left(\frac{\partial f}{\partial q}-\frac{\diff }{\diff t}\frac{\partial f}{\partial q'}+ \dots + (-1)^{k} \frac{\diff^{k} }{\diff t^{k}}\frac{\partial f}{\partial q^{(k)}}
  \right)\right|_{a(t)} , \delta _q(t) \right> {\rm d}t\\
  &+   &  \int_{I} \left<  \left. \left( \frac{\partial f}{\partial p}-\frac{\diff }{\diff t}\frac{\partial f}{\partial p'}+ \dots + (-1)^{l} \frac{\diff^{l} }{\diff t^{l}}\frac{\partial f}{\partial p^{(l)}} 
  \right)\right|_{a(t)}  , \delta _p(t) \right>  {\rm d}t\\
  &+   & \left< \left. \left( \frac{\partial f}{\partial q'}+ \dots + (-1)^{k} \frac{{\rm d}^{k-1}}{{\rm d}t^{k-1}} \frac{\partial f}{\partial q^{(k)}}\right)\right|_{a(1)}, \delta _q (1) \right> \\
  &-   & \left< \left. \left( \frac{\partial f}{\partial q'}+ \dots + (-1)^{k} \frac{{\rm d}^{k-1}}{{\rm d}t^{k-1}} \frac{\partial f}{\partial q^{(k)}}\right)\right|_{a(0)}, \delta _q (0) \right>.  
   \end{eqnarray*}
 To complete the previous computation, we simply used a successive integration by parts.
 By identification we obtain that $(A_F,B_F) $ exist and are given by:
 \begin{equation}
 \label{eq:AFBF}
 \left\{\begin{array}{rcl}
 A_F(t)&=& \left.\left( \frac{\partial f}{\partial q}-\frac{\diff }{\diff t}\frac{\partial f}{\partial q'}+ \dots + (-1)^{k} \frac{\diff^{k} }{\diff t^{k}}\frac{\partial f}{\partial q^{(k)}} \right)\right|_{a(t)}  \\
 B_F(t)&=& \left.\left( \frac{\partial f}{\partial p}-\frac{\diff }{\diff t}\frac{\partial f}{\partial p'}+ \dots + (-1)^{l} \frac{\diff^{l} }{\diff t^{l}}\frac{\partial f}{\partial p^{(l)}} \right)\right|_{a(t)} .
 \end{array}\right.
  \end{equation}
  \end{proof}
  
  \begin{remark}
 We saw in Remark \ref{rmk:non-unicite} that, given a local function $F$,  the function $f$ such that $F=F_f$ is not unique.
 However, the gradient at a given point of the local function $F$ is defined in a manner which is independent from the function $f$.
In particular, the above expressions for $A_F$ and $B_F$ do not depend on the choice of $f$. Notice that these expressions are linear in $f$.
  \end{remark}

  For loops, a similar proof yields the following statement:
\begin{lemma} \label{gradient-per}
For any local function $F \in \F({\mathfrak L})$, there exists one and only one element 
$\nabla _a F = (A_F,B_F ) \in \C^{\infty}(S^1,{\mathbb R}^{n}\times  {\mathbb R}^{n})$ satisfying that  for any tangent vector $\delta_a =(\delta_q,\delta_p)$ to ${\mathfrak L}$:
\begin{eqnarray*}\label{Lie bracket}
\left. \frac{\diff }{\diff \epsilon}\right|_{\epsilon=0} F(a+\epsilon \delta_a )
  &=   &\int_{I} \langle A_F ,\delta_q  \rangle {\rm d}t + \int_{I}   \langle B_F ,\delta_p  \rangle {\rm d}t. \\
   \end{eqnarray*}   
where $\langle . , . \rangle$ is the canonical scalar product of $ {\mathbb R}^n$.
Moreover $A_F$ and $B_F$ are given by (\ref{eq:AFBF}).
We call this pair the gradient of $F$ at $a$.
 \end{lemma}

  We define \cites{Ikeda,SchallerStrobl} a Lie bracket on local functions $F,G \in \F(\tilde{\mathfrak P})$ or $\F({\mathfrak L})$ as follows:
 \begin{equation}\label{eq:Liebracket}
\pb{F,G}(a)=  \int_{I} (\can{ A_F,B_G}-\can{ A_G,B_F} ){\rm d}t.
 \end{equation}
  where $ A_F,B_F,A_G,B_G$ are as in Lemmas  \ref{gradient-semi} and \ref{gradient-per} above. Although local functions do not form on algebra, we call this Lie bracket
  \emph{the canonical Poisson structure}. It is of course, up to a minor difference of setting, the Poisson bracket that appears in \cite{strobelklim,CattaneoFelder}.
  \begin{remark}
  The canonical Poisson structure is not symplectic on $\tilde{\mathfrak P}$: it admits Casimir functions. For instance, 
  for all smooth function $f$ on $M$
  the function $F(a) = \int_I \diff_{q(t)} f \, (q'(t)) \, {\diff t}$  (with $q$ the base path of $a$) is a local function on $\tilde{\mathfrak P}$ by definition.
  Since we have $F(a)= f(q(1)) - f(q(0))$, it satisfies $A_F=B_F=0$.
  Hence for all local function $G$, we have $\{F,G\}=0$, i.e. $F$ is a Casimir. 
  \end{remark}

We define cosiotropic subsets in an algebraic manner:   
 
 \begin{definition}\label{def:coisotropic}
  A subset $C$ of ${\mathfrak P}$ or ${\mathfrak L}$ shall be said to be coisotropic if and only if local functions vanishing on $C$ are stable under the canonical  Poisson structure 
  given in Equation (\ref{eq:Liebracket}).
  \end{definition}

\section{Non-periodic case}\label{Non-periodic case}
We introduce cotangent paths \cites{CrainicFernandes,FernandesCranicMarius}.
 \begin{definition}
Let $(M,\pi)$ be a manifold equipped with a bivector field $\pi$. 
A smooth path $a$ from $I=[0,1]$ to $T^*M$ is a \emph{cotangent path for $\pi$} if it satisfies:
  \begin{equation} \label{eq:def_chemins_cotangents}
 \pi_{q(t)}^\# (a) = \frac{\diff q(t)}{\diff t},
  \end{equation}
where $q$, called the \emph{base path of $a$}, is the projection of $a$ onto $M$.
We denote by $C_{\pi}$ the set of semi-free cotangent paths and 
by $C_{\pi}^{S_{1}}$ the set of loops which are cotangent paths.
 \end{definition}

\subsection{Theorem and proof}\label{Open sets}
 Again, we assume $M \subset {\mathbb R}^{n}$ is an open subset, so that paths $a$ in $ T^*M$
are pairs $q,p$ of paths in $M$ and ${\mathbb R}^n$ respectively, $q$ being the base path of $a$.  
 We describe a bivector field on ${\mathbb R}^{n}$ by a $n\times n$ skew-symmetric matrix:
$$
  \left(
     \raisebox{0.5\depth}{%
       \xymatrixcolsep{1ex}%
       \xymatrixrowsep{1ex}%
       \xymatrix{ 0 \ar@{.}[dddrrr] \ar@{.}[rrr] \ar@{.}[ddd]& & & \pi_{0n} \ar@{.}[ddd]\\ & & \pi_{ij}&  \\  & -\pi_{ij}& & \\- \pi_{0n} \ar@{.}[rrr]& & &  0 \\
       }%
     }
   \right)
  $$
whose coefficients are smooth functions from $M$ to ${\mathbb R}$ denoted by $\pi_{ij}$, for all $i,j \in \{1,\dots,n\}$.
For $t \mapsto q(t)$ a path in $M$, we denote in general by $ \pi_{ij}|_{q(t)}$ the value of this function at the point $q(t) \in M$. Of course,  the $s$-th component of $\pi^\#_{q(t)}(p(t)) $ is $ \sum_{j=1}^n \left. \pi_{sj}\right|_{q(t)}p_j(t)$.

 The set of cotangent paths can be characterized by local functions.
 \begin{remark}\label{rmk:caract_cotangent}
A path $a \in \tilde{{\mathfrak P}}$ is cotangent if and only if it is zero when evaluated on every local function of the form:
\begin{equation}\label{eq:coisotrope}
F_{f,s}(a)=\displaystyle{\int}_{I}  f(t) \left(\sum_{j=1}^n \left. \pi_{sj}\right|_{q(t)}p_j(t)-\frac{\diff q_{s}(t)}{\diff t}\right) \, \diff t .
\end{equation}
associated with the function on the space of $1$-jets given by $f_s(t,q,q',p,p') := f(t) \left(\sum_{j=1}^n \left. \pi_{sj}\right|_{q}p_j -q_{s} '\right)$, 
where $1 \leq s \leq n$ and where $f$ is a smooth function from $I$ to ${\mathbb R}$ vanishing at $t=0$ and $t=1$.
\end{remark}
 
We intend to prove the following result that we invite the reader to compare with the similar Equation (10) in \cite{SchallerStrobl}, Equation (7) in 
\cite{strobelklim}, Theorem 1.1 in \cite{Cattaneo} or Theorem 3.1 in \cite{CattaneoFelderCoisotropic}:
 
\begin{theorem}\label{theo:1}
Let $\pi$ be a bivector field on an open subset $M \subset {\mathbb R}^{n}$. The set $C_{\pi}$ of semi-free cotangent paths is a coisotropic subset of the set $\tilde{{\mathfrak P}}$ of semi-free paths if and only if $\pi$ is a Poisson bivector.
\end{theorem}

This result is a consequence of the two Propositions \ref{prop:si} and \ref{prop:ssi} that establish a sense then its reciprocal.
 
\begin{proposition}\label{prop:si}
If the set $C_{\pi}$ of cotangents paths is a coisotropic subset of the set $\tilde{{\mathfrak P}}$ 
of semi-free paths, then $\pi$ is a Poisson bivector.
\end{proposition}

We start by a lemma:
\begin{lemma}\label{lem:gradient}
Let $f$ a smooth function in $\C^{\infty}(I,{\mathbb R})$ and $1\leq s \leq n$. The gradient $\nabla F_{f,s}=(A_{F_{f,s}},B_{F_{f,s}})$ of the function $F_{f,s}$ defined in Remark \ref{rmk:caract_cotangent} 
at a point $a=(q,p) \in \tilde{\mathfrak P}$ is given by:
$$ A_{F_{f,s}} = f(t) \left( 
\begin{array}{c} \left.  \sum_{j=1}^n   \frac{\partial \pi_{sj}}{\partial q_1} \right|_{q(t)} \, p_j(t) \\ \vdots \\ \left.  \sum_{j=1}^n \frac{\partial \pi_{sj}}{\partial q_s} \right|_{q(t)} \, p_j(t)  \\ \vdots \\
\left. \sum_{j=1}^n  \frac{\partial  \pi_{sj}}{\partial q_n} \right|_{q(t)}  \, p_j(t)  \end{array}  \right) +  
f'(t) \left( \begin{array}{c}  0 \\  \vdots \\ 0 \\ 1 \\ 0 \\ \vdots \\ 0  \end{array} \right) \begin{array}{c}   \\   \\  \\ {\tiny{\leftarrow  (\hbox{s-th term)}}} \\  \\  \\ \\  \end{array}$$
and:
$$ B_{F_{f,s}}= f(t) \left( \begin{array}{c} \left. \pi_{s1} \right|_{q(t)} \\ \vdots \\  \left. \pi_{sn} \right|_{q(t)} \end{array} \right)  $$

\end{lemma}
\begin{proof}
The demonstration is given by a direct calculation. For all tangent vector to $\tilde{{\mathfrak P}}$, i.e. all elements $\delta _a$ in $\C^{\infty}(I,{\mathbb R}^{n}\times  {\mathbb R}^{n})$, we write:
$$\delta _a=(\delta _q,\delta _p)=\left(  \left(\begin{array}{c} \delta q_{1} (t)\\ \vdots \\ \delta q_{n} (t)\end{array} \right) ,
\left(\begin{array}{c} \delta p_{1} (t)\\  \vdots \\ \delta p_{n} (t) \end{array}  \right)\right)$$
for some functions $ \delta q_{1} ,\dots,\delta q_{n}, \delta p_{1} ,\dots,\delta p_{n}$ with appropriate boundary conditions. We have:
\begin{eqnarray*}
F_{f,s}(a+\epsilon \delta _a)
  &=   &\int_{I} f(t) \left(\sum_{j=1}^n \left. \pi_{sj}\right|_{q+\epsilon \delta _q} (p_j+\epsilon \delta p_j)-\frac{\diff q{_s}+\epsilon \delta _{q_s} }{\diff t}\right) \, \diff t  .
  \end{eqnarray*}
Differentiating the previous expression at $\epsilon=0$, and using integration by parts: 
   
\begin{eqnarray*}
\left.  \frac{\diff  F_{f,s}(a+\epsilon \delta _a)}{\diff \epsilon}\right|_{\epsilon=0} &=   & \int_{I}  f(t)  \can{ \left( 
\begin{array}{c}  \left. \pi_{sn}\right|_{q(t)}  \\  \vdots \\ \left. \pi_{sn}\right|_{q(t)} 
\end{array}  
\right) , \left(\begin{array}{c} \delta p_1 (t)\\  \vdots \\ \delta p_n (t) 
\end{array}  \right)  }  \, \diff t  \\
     &+   & \int_{I}  f(t) \can{ \left( 
\begin{array}{c} 
       \sum_{j=1}^n \left.\frac{\partial \pi_{sj}}{\partial q_1} \right|_{q(t)} p_j(t)
       \\ \vdots \\  
       \sum_{j=1}^n   \left.\frac{\partial \pi_{sj}}{\partial q_n} \right|_{q(t)} p_j(t) 
\end{array}    
  \right) , \left(\begin{array}{c} \delta q_1 (t)\\ \vdots \\ \delta q_n (t)
\end{array} \right)   }  \, \diff t \\
    &+   & \int_{I} f'(t) \delta q_s (t) {\rm d}t -  f(1) \delta q_s (1) + f(0) \delta q_s (0).  
\end{eqnarray*}
This gives the result by identification.
 \end{proof}
Let $\pi$ be a bivector field, whose associated biderivations on functions we denote by a $\PB$.
Recall that the map:
\begin{equation}
\begin{array}{rrcl}
Jac  :&C^{\infty}(M)^3  &\to & C^{\infty}(M)\\
      & (f,g,h)& \longmapsto &\pb{\pb{f,g},h}+\pb{\pb{h,f},g}+\pb{\pb{g,h},f},
\end{array}
\end{equation}%
is a $3$-vector field, called the \emph{Jacobiator} of $\pi$. It is zero if and only if $\pi$ is Poisson. We refer to \cite{ACP} for further details. In local coordinates, the Jacobiator is given by:
$$Jac(f,g,h) = \sum_{r,s,j=1}^n J_{rsj} \frac{\partial f}{\partial x_r} \frac{\partial g}{\partial x_s} \frac{\partial h}{\partial x_j}$$
with:
 \begin{equation}\label{eq:jacobiator} J_{rsj} = \sum_{k=1}^n \left( \frac{\partial \pi_{rs}}{\partial q_k} \pi_{kj} + \frac{\partial \pi_{sj}}{\partial q_k} \pi_{kr} +\frac{\partial \pi_{jr}}{\partial q_k} \pi_{ks} \right) .\end{equation}
 The following lemma follows from Lemma \ref{lem:gradient} (and is a particular case of Equation (7) in \cite{strobelklim}):
 \begin{lemma}\label{lem:crochet}
 Let $F_{f,r},G_{g,s}$ two functions of the form (\ref{eq:coisotrope}) where $f,g$ are two smooth functions vanishing at $t=0$ and $t=1$ in $C^{\infty}(I,{\mathbb R})$ and $r,s \in \{1,\dots,n \}$ . The bracket $\pb{F_{f,r},G_{g,s}}$ is given by:
 \begin{eqnarray*}\pb{F_{f,r},G_{g,s}}(a) &=& 
  \int_{I} f(t) g(t)  \sum_{k=1}^n \left. \frac{\partial \pi_{rs} }{\partial q_k}\right|_{q(t)} (q_k'(t)- \sum_{j=1}^n \pi_{kj}|_{q(t)}p_j(t)){\rm d}t 
 \\ &-& \int_{I} f(t) g(t) \sum_{j=1}^n J_{rsj}{|_{q(t)}} p_j(t){\rm d}t
 \end{eqnarray*}
at an arbitrary element  $a=(q,p)$ of $\tilde{{\mathfrak P}}$.
 \end{lemma}
 \begin{proof}
 By applying Lemma \ref{lem:gradient} and Equation (\ref{eq:Liebracket}):
\begin{eqnarray*}
& & \pb{F_{f,r},G_{g,s}}(q(t))\\
  &=   &\int_{I} (f'(t) g(t)  \pi_{sr}|_{ q(t)} - g'(t) f(t)  \pi_{rs}|_{q(t))}{\rm d}t\\
  &+   & \sum_{j,k=1}^n \int_{I} f(t) \left[ \left.\frac{\partial \pi_{rj}}{\partial q_k}\right|_{q(t)} p_j(t) g(t) \pi_{sk}|_{q(t)}-g(t) \left. \frac{\partial \pi_{sj}}{\partial q_k}\right|_{q(t)} p_j(t)  \pi_{rk}|_{q(t)} \right]{\rm d}t.
  \end{eqnarray*}
  The first of these two terms is, since $\pi_{sr}=-\pi_{rs}$:
  $$-\int_{I} (f'(t) g(t) +g'(t) f(t))  \pi_{rs}(q(t)) {\rm d}t.$$
  An integration by parts gives then:
  \begin{eqnarray*}
& & \pb{F_{f,r},G_{g,s}}(q(t))\\
  &=   &\int_{I} f(t) g(t)  \frac{\diff \pi_{sr}|_{q(t)} }{\diff t}{\rm d}t\\
  &+   & \sum_{j,k=1}^n \int_{I} f(t) \left[ \left.\frac{\partial \pi_{rj}}{\partial q_k}\right|_{q(t)} p_j(t) g(t) \pi_{sk}|_{q(t)}-g(t) \left.\frac{\partial \pi_{sj}}{\partial q_k}\right|_{q(t)} p_j(t)  \pi_{rk}|_{q(t)} \right]{\rm d}t  \\
  &=   &\int_{I} f(t) g(t)  \sum_{k=1}^n\frac{\partial \pi_{sr}}{\partial q_k} q_k'(t) {\rm d}t\\
  &+   & \sum_{j,k=1}^n \int_{I} f(t) g(t) p_j(t)  \left[ \left.\frac{\partial \pi_{rj}}{\partial q_k}\right|_{q(t)} \pi_{sk}|_{q(t)}- \left.\frac{\partial \pi_{sj}}{\partial q_k}\right|_{q(t)}  \pi_{rk}|_{q(t)} \right]{\rm d}t. 
  \end{eqnarray*}
Using Equation (\ref{eq:jacobiator}), this completes the proof of Lemma \ref{lem:crochet}. 

\end{proof}

The following lemma is an immediate consequence of Lemma \ref{lem:crochet}:
\begin{lemma}\label{lem:CCotang}
For all $r,s \in \{1, \dots,n \}$, and for all pair of smooth functions $ (f,g)$ vanishing at $t=0$ and $t=1$ and all semi-free cotangent path $a =(q,p) \in C_{\pi }$:
$$\pb{F_{f,r},F_{g,s}}(a) = \int_{I} f(t) g(t) \sum_{j=1}^n J_{rsj} p_j(t){\rm d}t ,$$
where $J$, the Jacobiator, is as in (\ref{eq:jacobiator}).
\end{lemma}

To conclude that the Jacobiator $J$ is zero, we need a technical lemma:

\begin{lemma}\label{lem:goesthrough}
 For all $(q,p)$ in $T^*M$, there exists a cotangent path $a $ such that  $a(1/2)=(q,p)$.
\end{lemma}
\begin{proof}
Let $(q,p)$ be an element of $ T^*M\simeq M \times {\mathbb R}^{n} $ where $M$ is open space of ${\mathbb R}^{n}$. We consider $y(t)$ a solution of the following differential equation defined for $t \in ]1/2-\epsilon ,1/2+\epsilon [$:
$$\pi^\#_{y(t)} (p) = \frac{\diff y(t)}{\diff t}~~\hbox{et}~~ y(1/2)=q.$$
By construction the map $t \mapsto  (y(t),p)$ is a cotangent path. We choose an application $\psi : I \to ]1/2-\epsilon ,1/2+\epsilon [$ such that $ \psi(1/2)=1/2$ and $\psi'(1/2)=1$ whose derivatives at $0$ and $1$ are equal to zero.
Then we consider the path $a : t \to (\psi '(t)p,y(\psi (t)))$. This path is semi-free by construction and it is cotangent path, since:
\begin{eqnarray*}
\frac{\diff }{\diff t}y(\psi (t)))
  &=   &\psi' (t) y'(\psi (t))\\
   &=   &\pi^\#_{y(t)}(\psi' (t)p).
 \end{eqnarray*}
By construction, $a(1/2)=(p,q)$. This shows the result.
\end{proof}
We can now prove Proposition \ref{prop:si}.

\begin{proof}
For all $(q,p) \in T^*M \simeq M \times {\mathbb R}^n$, we use Lemma \ref{lem:goesthrough} to choose a semi-free cotangent path $a$ such that $a (1/2)=(q,p)$.
Then we choose $g$ to be any smooth function vanishing at $0$ and $1$ and equal to $1$ in a neighborhood of $1/2$ and we chose $f=(f_d)_{d \geq 0}$ an approximation of the Dirac function at $t=1/2$. Then, we have:
 $$\lim_{d\rightarrow \infty}\int_{I} (f_d(t) g(t)) \sum_{j=1}^n J_{rsj}{|_{q(t)}} p_j(t){\rm d}t = \sum_{j=1}^n Jac(\pi )_{rsj}{|_{q(t)}} p_j$$
 where $p=(p_1, \dots,p_n)$.
By Lemma \ref{lem:CCotang}, we have therefore:
 $$\lim_{d\rightarrow \infty} \left(\pb{F_{f_d,r},G_{g,s}}(a) \right) =   \sum_{j=1}^n Jac(\pi )_{rsj}{|_{q(t)}} p_j .$$
If $C_\pi$ is coisotropic, then $\pb{F_{f,r},G_{g,s}}(a)=0$, so that
 $$\sum_{j=1}^n J_{rsj}{|_q} p_j=0.$$
As $q$ and $p$ are arbitrary, this gives that the Jacobiator $J$ vanishes, i.e. that $\pi$ is Poisson. 
\end{proof}

We now prove the reverse implication.

\begin{proposition}\label{prop:ssi}
If $\pi$ is a Poisson bivector, then $C_{\pi}$ is coisotropic.
\end{proposition}
We start by a definition.

\begin{definition}
Let $a=(q,p) \in T^{*}M\simeq M\times {\mathbb R}^{n}$ a semi-free cotangent path.
We say that the vector $(\delta_q ,\delta_p ) \in C^{\infty}(I,{\mathbb R}^{2n})$ is tangent to $C_{\pi}$ at the point $a$ if there exists a path $a_{\epsilon }=(q_{\epsilon },p_{\epsilon })$ valued in $C_{\pi}$ where $q_{\epsilon }(t)$ and $p_{\epsilon }(t)$ are smooth functions in $t$ and $\epsilon $ such that:
\begin{enumerate}
\item[1)] $q_{0}=q$ and $p_{0}=p$,
\item[2)] $\frac{\diff }{\diff \epsilon }_{|_{\epsilon =0}}q_{\epsilon }=\delta _q$ and $\frac{\diff }{\diff \epsilon }_{|_{\epsilon =0}}p_{\epsilon }=\delta _p$
 We denote by $TC_{\pi}$ and call tangent cone of $C^{\pi}$ at $a$ the set of elements tangent to $C_\pi$. 
\end{enumerate}
\end{definition}
A priori it is not obvious that $TC^{\pi}$ is a vector space. This is true, however, in view of the following lemma.
Let us introduce a notation. For all $u,p \in {\mathbb R}^n$, we denote by:
 $$\frac{\partial{\pi^\# }}{{\partial q}}(u)(p) $$
 the derivative of $q \mapsto \pi^\#|_q (p)$ in the direction of $u$.
 Equivalently, it is the derivative in the direction of $u$ of the $n\times n$-matrix $\pi^\#$, to which we apply $ p \in {\mathbb R}^n$.

\begin{lemma}\label{lem:tangentCpi}
Let $a=(q,p) \in C_{\pi}$ a semi-free cotangent path and $\delta_a= (\delta_q ,\delta_p ) \in C^{\infty}(I,{\mathbb R}^{n} \times {\mathbb R}^n)$ a vector tangent to ${\tilde {\mathfrak P}}$ at the point $a \in C_{\pi}$.
We have  $\delta_a \in TC_{\pi}$ if and only if $\frac{\diff}{\diff t}(\delta _q)=\frac{\delta \pi^\#}{\delta q}(\delta _q)(p)+ \pi^\#_q (\delta _p).$
In particular $TC_{\pi}$ is a vector space.
\end{lemma}

\begin{proof} 
Let $a=(q,p) \in C_{\pi}$ a semi-free cotangent path and $(\delta_q ,\delta_p ) \in C^{\infty}(I,{\mathbb R}^{2n})$ a vector tangent to $C_{\pi}$ at the point $a$.
We consider $p_{\epsilon} = p + \epsilon \delta_p$ and $q_{\epsilon}(t)$ the flow of the differential equation:
$$q_{\epsilon}'=\pi_{q_{\epsilon}(t)} (p_{\epsilon})~~~~~~ \hbox{with}~~~~~ q_{\epsilon}(0) = y(\epsilon) ,$$
where $y(\epsilon)$ is as $y(0) = q_{0}$ and $y'(0) = \delta_q(0)$.
For all $\epsilon $, the path $ t \mapsto (q_{\epsilon}(t), p_{\epsilon}(t))$ is a cotangent path in $C_{\pi}$ by construction. 
Also, since $ p + \epsilon \delta_p$ vanishes together with all its derivatives at $t=0$ and $t=1$, 
so does $\frac{\diff}{\diff t} q_\epsilon$, so, for all $\epsilon $, the path $ t \mapsto (q_{\epsilon}(t), p_{\epsilon}(t))$ is semi-free. 
Also, by construction again, we have $(q_{0}, p_{0})=(q,p)=a$. Hence 
$\frac{\diff }{\diff \epsilon }(q_{\epsilon}, p_{\epsilon})_{|_{\epsilon =0}} = \left(\left.\frac{\diff q_{\epsilon}}{\diff \epsilon}\right|_{\epsilon=0},\delta_p \right)$ is tangent to $C_{\pi}$ at $a$.
Since we have by construction:
$$\frac{\diff q_{\epsilon}}{\diff t }=\pi_{q_{\epsilon}}^\# (p_{\epsilon}),$$
 it turns that
$$\frac{\diff}{\diff t} \left( \left. \frac{\diff q_{\epsilon}}{\diff \epsilon } \right|_{\epsilon =0} \right)=
\frac{\delta \pi^\#}{\delta q} \left( \left.\frac{\diff q_{\epsilon}}{\diff \epsilon }\right|_{\epsilon =0}\right)(p)+
\pi^\#_{q} \left( \left.\frac{\diff p_{\epsilon } }{\diff \epsilon } \right|_{\epsilon =0}\right)+\pi^\#_q (\delta _p).$$
Since $\delta _q$ and $\left.\frac{\diff q_{\epsilon}}{\diff \epsilon } \right|_{\epsilon =0}$ satisfy the same differential equation with the same initial value, they are equal.
In conclusion,  $\delta_a=(\delta_q, \delta_p) \in T C_\pi$. 
\end{proof}

Let $\varphi$ be the assignment given by:
\begin{equation}\label{eq:diff}
\renewcommand{\arraystretch}{1.7}
\begin{array}{cccc}
\varphi   :&C^{\infty}(I,{\mathbb R}^{n} \times {\mathbb R}^{n})&\to& C^{\infty}(I,{\mathbb R}^{n})\\
      &(q,p)&\mapsto& q'(t)-\pi^\#_{q(t)} (p).
\end{array}
\end{equation}%
Of course, $C_{\pi}$ is the restriction to semi-free paths of $\varphi ^{-1}(0)$. 
The differential of $\varphi  $ is given by:
\begin{eqnarray*}
\diff \varphi   (\delta _q,\delta _p)
  &=   & \delta _q'-\frac{\delta  \pi}{\delta q} ( \delta _q) (p)+\pi (\delta _p).
\end{eqnarray*}
\begin{remark}
By Lemma \ref{lem:tangentCpi}, for all $a \in C_\pi$, the tangent cone $T_a C_{\pi} $ is the kernel $ Ker (\diff_a \varphi )$ of the differential of $\varphi$ at $a$.
\end{remark}
The application $\varphi$ of (\ref{eq:diff}) is an application that can be defined for $(q,p)$ in $H^{1}(I,{\mathbb R}^{n} \times {\mathbb R}^n)$.
It is then valued in $L^{2}(I,{\mathbb R}^{n})$.
In fact, if $q \in H^{1}(I,{\mathbb R}^{n} )$ then $q$ is bounded and consequently $\pi_{q(t)}^\#$ is bounded independently of $t \in I$. Since $p \in H^{1}(I,{\mathbb R}^{n} )$, we also have
that $p \in L^{2}(I,{\mathbb R}^{n})$ and consequently $\pi^\#(p)$ is also in $L^{2}(I,{\mathbb R}^{n})$.
  \begin{remark}\label{rem:existance}
Let $F\in\F (\tilde{{\mathfrak P}})$ a local function vanishing on $C_{\pi}$ and $a \in C_\pi$. We know that $\diff_a F(\delta_a)$ is zero if $\delta_a \in TC_{\pi}$, in particular $\diff F(\delta_a)=0$ if $\delta_a \in \Ker \diff_a {\varphi}$. 
 By the Hahn-Banach theorem, there is an element $H \in  L^2(I,{\mathbb R}^n)$, such that the following diagram is commutative:
\begin{equation*}\label{diag:bana}
 \begin{diagram}
\node{H^{1}(I,{\mathbb R}^{n} \times {\mathbb R}^n)}\arrow[2]{e,t}{\diff_a \phi}\arrow{se,b}{\diff_a F}\node[2]{L^{2}(I,{\mathbb R}^{n})}\arrow{sw,b}{\can{H,.}}\\
\node[2]{{\mathbb R}}
\end{diagram}
\end{equation*}%
where $\can{.,.}$ is the natural scalar product of $L^{2}(I,{\mathbb R}^{n})$.
\end{remark}
From Remark \ref{rem:existance}, we deduce the following lemma:
\begin{lemma}\label{lem:Hbanach}
Let $F \in \F (\tilde{{\mathfrak P}})$ a local function vanishing in $C_\pi$, let $a=(q,p)$ a cotangent path and $\delta_a=(\delta_q ,\delta_p )$ a tangent vector to $\tilde{\mathfrak P}$.
Then there exists a smooth function $H:I \to {\mathbb R}^{n}$ such that the components $A_F,B_F$ of the gradient of $F$ at $a$ are given by:
\begin{equation}\label{equa:exisH}
\left\{
\begin{array}{rcl}
A_F(t) &= & -H'(t)+ \frac{\partial \pi}{\partial q}^*|_{q(t)} (p(t),H(t))\\
B_F(t) & =& \pi^\#|_{q(t)} H(t)\\
\end{array}
 \right.
\end{equation}
where $\frac{\partial \pi}{\partial q}^*$ is at a given $ q \in M$ the bilinear assignment on ${\mathbb R}^n$ given by:
 $$  \left\langle \left.\frac{\partial \pi}{\partial q}^*\right|_q (\alpha,\beta) , u \right\rangle = -\left\langle \left.\frac{\partial \pi^\#}{\partial q} \right|_{q} (u)(\beta) ,\alpha  \right\rangle,$$
for all $\alpha,\beta ,u \in {\mathbb R}^n$. 
\end{lemma}
\begin{proof}
From the diagram (\ref{diag:bana}) we know that there exists a function $H \in L^{2}(I,{\mathbb R}^n)$ such that for all $\delta_a=(\delta_q ,\delta_p )$ tangent to $\tilde{\mathfrak P}$ at $a$ with $ \delta_q=0$ in neighborhoods of $0$ and $1$, we have:
  \begin{eqnarray*}
\diff_a F(\delta _q,\delta _p)
  &=   & \int_{I} \can{H,\diff_a \phi  (\delta _q,\delta _p}) {\rm d}t\\
   &=   &\int_{I} \can{H,\delta _q'- \left. \frac{\delta  \pi^\#}{\delta q}\right|_{q(t)} (\delta _q) (p) -\pi^\# \delta _p} {\rm d}t\\
    &=   &\int_{I} \can{H, \delta _q'- \left. \frac{\delta  \pi^\#}{\delta q}\right|_{q(t)} (\delta _q) (p)} {\rm d}t+ \int_{I} \can{\pi^\# H, \delta _p} {\rm d}t\\
      &=   &  \int_{I} \can{-H'+\frac{\delta  \pi}{\delta q}^* (p, H),\delta _q }+ \int_{I} \can{\pi^\# (H),\delta _p} {\rm d}t.
  \end{eqnarray*}
Here $H'$ is the derivative of the function $H$ in the sense of distributions.
By definition of the gradient, we find an identification (in the sense of distributions) which is precisely Equation (\ref{equa:exisH}).
Since $A_F$ is a smooth function, and $H$ is in $L^{2}(I,{\mathbb R}^n)$, $H'$ is in $L^{2}(I,{\mathbb R}^n)$ and therefore $H$ is in $H^{1}(I,{\mathbb R}^n)$.
In particular, it is continuous which implies that $H$ is differentiable  and, by an immediate recursion, $H$ is a smooth function. By construction,
 Equation (\ref{equa:exisH}) holds.
\end{proof}

Last, we need to reinterpret Lemma \ref{lem:gradient2}:

\begin{lemma}\label{lem:gradient2}
Let $f$ a smooth function in $\C^{\infty}(I,{\mathbb R})$ and $1\leq s \leq n$. The gradient $\nabla F_{f,s}=(A_{F_{f,s}},B_{F_{f,s}})$ of the function $F_{f,s}$ defined in Remark \ref{rmk:caract_cotangent} 
at a point $a=(q,p) \in \tilde{\mathfrak P}$ is given by:
$$ A_{F_{f,s}} = f_s'(t)- \left. \frac{\partial \pi}{\partial q}^*\right|_{q(t)} (p(t),f_s(t)) $$
and
$$ B_{F_{f,s}}=  \pi^\# f_s $$
where
$$ f_s=  \left( \begin{array}{c}  0 \\  \vdots \\ 0 \\ f(t) \\ 0 \\ \vdots \\ 0  \end{array} \right) \begin{array}{c}   \\   \\  \\ {\tiny{\leftarrow  (\hbox{s-th term)}}} \\  \\  \\ \\  \end{array} $$
\end{lemma}

Now we prove the Proposition \ref{prop:ssi}.

\begin{proof}
Let $F$ a local function vanishing in $C_{\pi}$ and $a=(q,p) \in C_\pi$.
By Lemma \ref{lem:Hbanach} there exists a smooth function $H=(h_{1},\dots,h_{n})$ in $C^\infty(I,{\mathbb R}^n)$ such that: 
(\ref{equa:exisH}) holds true. By comparing Lemma \ref{lem:gradient} and Lemma \ref{lem:gradient2},
we see that $F$ and $-\sum_{i=1}^n F_{h_i,i} $ have the same gradient at the point $a$.

Hence, for $F,G$  a pair of local functions vanishing on $C_{\pi}$ and for all $a \in C_\pi$,
 there exists a family of smooth functions $(F_i,G_i)$ where $i$ is a in index between $1$ and $n$ such that:
$$ \begin{array}{rr} A_{F}= -\sum_{i=1}^n A_{F_{f_i,i}} &  B_{F}=- \sum_{i=1}^n B_{F_{f_i,i}} \\ A_{G}= -\sum_{j=1}^n A_{G_{g_j,j}} & B_{G}=- \sum_{j=1}^n B_{G_{g_j,j}}. 
\end{array}
 $$
Consequently, in view of Equation (\ref{eq:Liebracket}):
\begin{eqnarray*}
\pb{F,G}(a)
  &=   & \sum_{i,j=1}^n \pb{F_{f_i,i},G_{g_j,j}} (a).
 \end{eqnarray*}
By Lemma \ref{lem:crochet}, if $\pi$ is Poisson, then this Poisson bracket is zero evaluated at all cotangent semi-free paths. This implies that $\pb{F,G}(a)=0$ and completes the proof.
\end{proof}

\section{Periodic case}\label{periodiccase}

\subsection{Theorem and proof}\label{lecasperiodique}

From now on, we turn our attention to the case of smooth loops valued in the cotangent space of a manifold $M$. 
Recall from Section  \ref{sec:local} that smooth loops are denoted by ${\mathfrak L}$ and that local functions on ${\mathfrak L}$ are restriction to ${\mathfrak L}$ of local functions on 
${\mathfrak P}$. 
Recall from Section \ref{sec:local} that $C_{\pi}^{S^1}$ stands for loops which are cotangent paths, that are simply referred as cotangent loops.
\begin{theorem}\label{theo:2}
Let $\pi$ be a bivector field in $M \subset {\mathbb R}^{n}$. The set $C_{\pi}^{S^1}$ of cotangent loops is a coisotropic subset of the set ${\mathfrak L}$ of
 smooth loops valued in $T^*M$ if and only if $\pi$ is a Poisson bivector.
\end{theorem}

We now devote the rest of this section to a proof of this result. As in the semi-free case, we first turn our attention to the case where $M $ is an open subset of ${\mathbb R}^d$.

\begin{remark}
Theorem \ref{theo:2} is inspired by Theorem 1.5 in \cite{Cattaneo}, where a similar result is stated for loops of class $C^1$.
The statement seems to us to be slightly incorrect. More precisely,  Theorem 1.5 in \cite{Cattaneo} states that the set $C_{\pi}^{S^1}$ of cotangent loops is a coisotropic submanifold of ${\mathfrak L}$ of loops 
of class $C^1$ valued in $T^*M$ if and only if $\pi$ is a Poisson bivector. Unfortunately, $ C_{\pi}^{S^1}$ may not be a Banach submanifold (see Section \ref{sec:counterExample}) 
so that we can only speak of a 
coisotropic subset. 
\end{remark}

One sense of the proof of theorem \ref{theo:2} goes exactly as in the non-periodic case: if the set $C_{\pi}^{S^1}$ of cotangent loops is a coisotropic subset of the set ${\mathfrak L}$ of
 smooth loops valued in $T^*M$, then $\pi$ is a Poisson bivector. The proof of the semi-free case (i.e. Proposition \ref{prop:si}) can be adapted step by step and, using the functions introduced
 in Remark \ref{rmk:caract_cotangent}, we arrive at a conclusion similar to Lemma \ref{lem:CCotang}:
 
 \begin{lemma}\label{lem:CCotang2}
For all $r,s \in \{1, \dots,n \}$, and for all pair of periodic smooth functions $ (f,g) \in {\C}^\infty(S^1,{\mathbb R})$ and all cotangent loop $a =(q,p) \in C_{\pi }^{S^1}$:
$$\pb{F_{f,r},F_{g,s}}(a) = \int_{I} (f(t) g(t)) \sum_{j=1}^n J(\pi )_{rsj} p_j(t){\rm d}t.$$
\end{lemma}
\noindent
But Lemma \ref{lem:goesthrough} has to be adapted:

\begin{lemma}
 For all $(q,p)$ in $T^*M$, there exists a cotangent loop $a $ such that  $a(1/2)=(q,p)$.
\end{lemma}
\begin{proof}
Let $(q,p)$ a element of $ T^*M\simeq \U \times {\mathbb R}^{n} $ where $\U$ is open space of ${\mathbb R}^{n}$. We consider $y(t)$ a solution of the following differential equation defined for $t \in ]1/2-\epsilon ,1/2+\epsilon [$ :
$$\pi^\#_{y(t)} (p) = \frac{\diff y(t)}{\diff t}~~\hbox{et}~~ y(1/2)=q.$$
By construction the maps $t \mapsto  (y(t),p)$ is a cotangent path. Define an application $\psi : S^1 \to ]1/2-\epsilon ,1/2+\epsilon [$ such that $\psi '(1/2)=1$
and $\psi (1/2)=1/2$.
Then consider  $a : t \to (\psi  '(t)p,y(\psi  (t)))$. This is a periodic cotangent path, since:
\begin{eqnarray*}
\frac{\diff y(\varphi (t))}{\diff t}
  &=   &\varphi' (t) y'(\varphi (t))\\
   &=   &\pi^\#_{y(t)}(\varphi' (t)p).
 \end{eqnarray*}
and by construction $a(1/2)=(p,q)$. This shows the result.
\end{proof}

The rest of the proof of the "if" direction in Theorem \ref{theo:2} goes exactly as in the semi-free case. Let us now establish the other direction, similar to
Proposition \ref{prop:ssi}:

\begin{proposition}\label{prop:ssi2}
Let $\pi$ be a bivector field in $M \subset {\mathbb R}^{n}$. If if $\pi$ is a Poisson bivector, then the set $C_{\pi}^{S^1}$ of cotangent loops is a coisotropic subset of the set ${\mathfrak L}$ of
 smooth loops valued in $T^*M$.
\end{proposition}

The difficulty is that, as we shall see in Section \ref{sec:counterExample}, cotangent loops do not for a submanifold in general, so that it is complicated to characterize local functions vanishing on it, 
and their gradients. The idea of the proof, however, is that cotangent loops are a disjoint union, indexed by symplectic leaves of $\pi$, of cotangent loops 
over a given symplectic leaf. Now, we shall prove that for every symplectic leaf, the cotangent loops over that leaf form (an infinite dimensional equivalent of) a Lagrangian submanifold and is a coisotropic subset. 
Let us be more precise. For $S$ a symplectic leaf of $\pi$, denote by $C_\pi^{S^1}(S)$ cotangent loops whose base paths are valued in $S$.
By construction, the base path of a cotangent path is valued in some symplectic leaf of $\pi$,
it can not "jump" from a symplectic leaf to an other one; hence:
 $$ C_\pi^{S^1} = \coprod_{S \in {\mathcal S}} C_\pi^{S^1}(S), $$
 with $ {\mathcal S}$ the set of symplectic leaves of $\pi$. Since an arbitrary union of coisotropic subsets is a coisotropic subset, Proposition \ref{prop:ssi2}, and therefore Theorem \ref{theo:2}, 
 follows directly from the following result.

\begin{proposition}\label{prop:ssi3}
Let $\pi$ be a Poisson bivector field in $M \subset {\mathbb R}^{n}$ and $S$ a symplectic leaf. Then $ C_\pi^{S^1}(S)$ is a coisotropic subset.
\end{proposition}

From now on, we fix a symplectic leaf $S$, and we fix some sub-bundle ${ \mathfrak j}$ of $TM|_S$ in direct sum with $TS$:
 $$ TM|_S={ \mathfrak j} \oplus TS ,$$
then we decompose $T^*M|_S$ as a direct sum:
 $$ T^*M|_S = {\mathfrak h} \oplus {\mathfrak k} $$
with $ {\mathfrak k}:={\rm Ker}(\pi^\#) $ and $ {\mathfrak h}={\mathfrak j}^{\perp} $ (i.e. covectors vanishing on ${\mathfrak j}$)).
Also, we use $\langle . , . \rangle$ the canonical scalar product of $ {\mathbb R}^n$ to identify tangent and cotangent spaces of $M$ to ${\mathbb R}^n$.
Of course, the restriction of $\pi^\#$ to the symplectic leaf $S$ is of constant rank, which justifies the definition of ${\mathfrak k}$. 
By construction, $\pi^\#$ restricts to a vector bundle isomorphism between ${\mathfrak h}$ and ${TS}$.
For any cotangent loop $a=(p,q)$ over the symplectic leaf $S$, the tangent space of ${\mathfrak L}$ at the point $a$ can therefore by identified as quadruples
$(u,v,k,h)$ making the following diagram commutative:
\begin{equation}\label{eq:quadruples}
\xymatrix{ S^1 \ar[rrd]^{q} \ar[rr]^{(u,v,k,h)} & & \ar[d] TS \oplus {\mathfrak j} \oplus {\mathfrak k} \oplus {\mathfrak h} \\ & & S}
 \end{equation}
i.e. quadruples $(u,v,k,h)$ in $C^\infty (S^1,TS) \oplus C^\infty (S^1,{\mathfrak j} )\oplus C^\infty (S^1,{\mathfrak k}) \oplus C^\infty (S^1,{\mathfrak h})$ 
over the base path $q: S^1 \to S$. Let us relate this identification with the obvious identification of an element $\delta_a \in T_a {\mathfrak L} $ as a pair
$ \delta_a = (\delta_q, \delta_p) \in C^\infty(S^1,{\mathbb R}^n ) \times C^\infty (S^1 ,{\mathbb R}^n )$, which arises from the identification $T^*M =
M \times {\mathbb R}^n \subset {\mathbb R}^n  \times {\mathbb R}^n $. It is clear that 
\begin{equation}\label{identificationTangentSpace} \delta_q = u+v \hbox{ and }  \delta_p = h+k \end{equation}
 From now on, elements in $T_a {\mathfrak L}$
 shall simply be denoted as quadruples of the form  $(u,v,k,h) $ with the conventions above.
  To start with, we intend to prove the following result:

\begin{lemma}\label{lem:tgtconeIsSpace}
For all $a=(q,p) \in C_\pi(S)$ tangent cone of $C_\pi^{S^1}(S)$ at $a$ is the subspace of $T_{a}{\mathfrak L}$  made of all $(u,v,k,h) \in T_a {\mathfrak L}$
with $v=0$ and
 \begin{equation}\label{eq:tangentSpace}
 \pi^{\#} h  = u' -  \frac{\partial \pi^{\#}}{\partial q}(u)(p).
  \end{equation}
  In particular, for every maps $u$ and $k $ from $S^1$ to respectively $TS$ and ${\mathfrak k}$ over $q:S^1 \to S$, there exists an element in $T_a C_\pi(S)$ of the form $(u,0,k,h)$
\end{lemma}
Notice that neither $u'$ nor $\frac{\partial \pi^{\#}}{\partial q}(u)(p)$ are a priori sections of $TS$, but their difference has to be a section of $TS$.

\begin{proof}
For every $a=(q,p) \in C_\pi(S)$ and every one-parameter family $ a_\epsilon=(q_\epsilon,p_\epsilon)$ of elements in $C_\pi(S)$ with $a_0=a$, the derivative at $\epsilon=0$ of $ q_\epsilon$ is a loop valued in $TS$ (over the loop $q$). This already proves that the derivative at $\epsilon=0$ of  $a_\epsilon$ is of the form $(k,h,u,0)$.
Now, by differentiating at $\epsilon=0$ the relation:
 $$  q_{\epsilon}'= \pi^{\#}_{q_\epsilon} p_{\epsilon} $$
we obtain Condition (\ref{eq:tangentSpace}). This gives one part of the proof. 
We now prove the converse. Consider an element in $T_a{\mathfrak L}$ of the form $(u,v,k,h)$ with $v=0$ and where $h$ and $u$ satisfy Condition (\ref{eq:tangentSpace}). 
Since $v=0$, and since $S$ is a manifold, the exists a one-parameter family $q_{\epsilon}$ of loops valued in $S$ with $q_0=q$ whose derivative with respect 
to $\epsilon$ at $\epsilon=0$ is $u \in C^\infty_q(S^1,TS)$. Since $\pi^{\#}: {\mathfrak h} \to TS$ is one-to-one, there exists, for every $\epsilon$, an unique loop $p_\epsilon$ valued in
$ {\mathfrak h}$ such that $ \pi_{q_\epsilon}^{\#} p_\epsilon = q_\epsilon'$. The henceforth obtained family $a_\epsilon=(q_\epsilon,p_\epsilon)$
is one-parameter family valued in $C_\pi(S)$. Of course, adding to $p_\epsilon$ an arbitrary loop $k$ valued ${\mathfrak k}_\epsilon$, we still obtain a one-parameter family valued in $C_\pi(S)$ that we still denote by $a_\epsilon$. Moreover, upon choosing $k$ such  that $k_0+h_0=p$ (which is always possible since $TM|_S= {\mathfrak k} \oplus {\mathfrak h}$),the henceforth obtained family of loops $a_\epsilon$ goes through $a=(q,p)$ at $\epsilon=0$
and is valued in $C_\pi(S)$. According to first part of the proof, the derivative at $\epsilon$ of this family of loops is of the form:
 $(u,0,\tilde{k},\tilde{h})$ where $\tilde{h}$ and $u$  satisfy Condition (\ref{eq:tangentSpace}), which implies that $\tilde{h}=h$,
 since $h$ and $u$ also satisfy Condition (\ref{eq:tangentSpace}) and $\pi^{\#}$ is one-to-one. In general, we do not have $\tilde{k}=k$, but, 
 since we can add to $a_\epsilon $ any one-parameter family of loops $k_\epsilon$ valued in ${\mathfrak k}$ provided that $k_0=0$ 
 and still get one-parameter family of elements in $C_\pi(S)$ through $a=(q,p)$, we can obtain such a family whose derivative at $\epsilon=0$ is $(u,v,k,h)$.
\end{proof}

There is a natural skew-symmetric bilinear assignment on $T{\mathfrak L} $ which, under identifying the identification of an element $\delta_a \in T{\mathfrak L}$ 
as a pair  $ ( \delta_q, \delta_p ) \in C^\infty(S^1,{\mathbb R}^n) \times  C^\infty(S^1,{\mathbb R}^n)$ is given,
for all $\delta_{a_1} =  ( \delta_{q_1}, \delta_{p_1} )$ and $\delta_{a_2} =  ( \delta_{q_2}, \delta_{p_2} )$, by:
 $$ \omega ( \delta_{a_1},\delta_{a_2}) = \int_I (\left\langle \delta_{p_1} , \delta_{q_2} \right\rangle - \left\langle  \delta_{p_2} , \delta_{q_1} \right\rangle ) {\rm d}t ,$$
 where $\langle . , . \rangle$ is the canonical scalar product of $ {\mathbb R}^n$. The next lemma gives the expression
of the skew-symmetric bilinear assignment $\omega$ when $\delta_{a_1} $ and $\delta_{a_2} $ are expressed as quadruples  $(u_{1},v_{1},k_{1},h_{1})$ and $(u_{2},v_{2},k_{2},h_{2})$,
as in (\ref{eq:quadruples}), respectively.
  
\begin{lemma}\label{lem:scalar}
For all $a=(q,p) \in {\mathfrak L}$, for every two $\delta_{a_1} = (u_{1},v_{1},k_{1},h_{1})$ and $\delta_{a_2}  = (u_{2},v_{2},k_{2},h_{2})$ in $ T_a {\mathfrak L}$,
we have : 
$$\omega (\delta_{a_1}, \delta_{a_2}) = \int_{S^1}  \left( \langle k_{1} ,v_{2}  \rangle - \langle k_{2} ,v_{1} \rangle + \langle h_{1} ,u_{2}   \rangle - \langle h_{2} ,u_{1}   \rangle \right) {\rm d}t.$$
\end{lemma}
\begin{proof}
The proof simply follows from the fact that ${\mathfrak j}$ and ${\mathfrak h}$ are orthogonal one to the other by construction, and from the skew-symmetry of $\pi^\#_m : T^*_m M \to T_m M$ at all $m \in M$, which implies that its image (which is $TS$) is orthogonal to its kernel (which is ${\mathfrak k}$). 
\end{proof}

 Also, we need a characterization of the Jacobi identity:
 
\begin{lemma}\label{lemma:Jacobi}
A bivector field $\pi$ on $M \subset {\mathbb R}^n$
satisfies the Jacobi identity if and only if for all covectors $ \alpha_1,\alpha_2,p \in T^*_m M \simeq {\mathbb R}^n$ at some point $m \in M$,
 we have:
 $$
\left\langle \frac{\partial \pi^\#}{{\partial q}} (\pi^\# p) (\alpha_2) - \frac{\partial \pi^\#}{{\partial q}} (\pi^\#(\alpha_2)) (p) , 
\alpha_{1}   \right\rangle = - \left\langle   \frac{\partial \pi^\#}{{\partial q}} (\pi^\#  \alpha_1) (p), \alpha_{2}   \right\rangle $$
\end{lemma}
\begin{proof}
It follows from Equation (\ref{eq:jacobiator}) that the Jacobi identity is equivalent to:
$$  \left\langle \frac{\partial{\pi^\#}}{{\partial q}} (\pi^\# p ) (\alpha_1) , \alpha_2   \right\rangle  
 + \left\langle \frac{\partial{\pi^\#}}{{\partial q}} (\pi^\# \alpha_1 ) (\alpha_2) , p   \right\rangle 
 + \left\langle \frac{\partial{\pi^\#}}{{\partial q}} (\pi^\# \alpha_2 ) (p) , \alpha_1   \right\rangle 
  =0. $$
Lemma \ref{lemma:Jacobi} then follows by skew-symmetry of $\pi^\#$.
\end{proof}

We can now prove the following lemma, which should be interpreted as meaning that $T_a C_\pi^{S^1} (S)$ is Lagrangian for $\omega$, i.e.
that $C_\pi^{S^1} (S) $ is a Lagrangian submanifold.

\begin{lemma}
For all $a=(q,p) \in C_\pi^ {S^1} (S)$ a  cotangent path over $S$. 
\begin{enumerate}
\item For all  $\delta_1 a , \delta_2 a \in  C_\pi^{S^1} (S)$, we have $\omega ( \delta_1 a,\delta_2 a ) = 0$.
\item Conversely, every  $\delta_1 a \in  T_a {\mathfrak L}$ such that $\omega ( \delta_1 a,\delta_2 a ) = 0$ for all $\delta_2 a \in T_a  C_\pi^{S^1} 
(S)$ is an element in $T_a  C_\pi^{S^1} (S)$.
\end{enumerate}
\end{lemma}
\begin{proof}
 For any two $\delta_{a_1} = (k_1,h_{1},u_{1},v_1)$ in  $T_a^{S^1}{\mathfrak L}$ and $\delta_{a_2} = (k_2,h_{2},u_{2},v_2)$ in $C_\pi^{S^1} (S)$  (which
  means $v_2=0$ while $h_{2},u_{2}$ satisfies Condition (\ref{eq:tangentSpace}) by Lemma \ref{lem:tgtconeIsSpace}), we have in  view of Lemma \ref{lem:scalar}:
\begin{eqnarray*}\omega (\delta_{a_1}, \delta_{a_2}) &=& \int_{S^1}  \left(\langle h_{1} ,u_{2}   \rangle - \langle h_{2} ,u_{1}   \rangle -
   \langle k_2, v_1 \rangle
   \right) {\rm d}t \\ 
     & = & \int_{S^1}  \left(\left\langle h_{1} ,\pi^\# (\alpha_{2})   \right\rangle - \left\langle h_{2} , \pi^\# (\alpha_{1})   \right\rangle 
   - \langle k_2, v_1 \rangle \right) {\rm d}t \\
   & = & -\int_{S^1}  \left(\left\langle \pi^\# ( h_{1} ) ,\alpha_{2}   \right\rangle - \left\langle \pi^\# (h_{2}) , \alpha_{1}   \right\rangle \right) {\rm d}t
   -\int_{S^1}  \langle k_2, v_1 \rangle  {\rm d}t \\
& = &  -\int_{S^1}  \left( \left\langle \pi^\# (h_{1}) , \alpha_{2}   \right\rangle - 
\left\langle \frac{\diff u_2}{\diff  t}- \frac{\partial {\pi^\#}}{{\partial q}} (u_2) (p) , \alpha_{1}   \right\rangle  \right) {\rm d}t \\
&-& \int_{S^1} \langle k_2, v_1 \rangle  {\rm d}t \\
& = &  -\int_{S^1} \left\langle \pi^\# (h_{1}) , \alpha_{2}   \right\rangle {\rm d}t \\
 & + & \int_{S^1}  \left\langle \frac{\diff \pi^\#(\alpha_2)}{\diff  t}- \frac{\partial{\pi^\#}}{{\partial q}} (\pi^\#(\alpha_2)) (p) , \alpha_{1}   \right\rangle   {\rm d}t- \int_{S^1} \langle k_2, v_1 \rangle  {\rm d}t,
  \end{eqnarray*}
where $\alpha_i,i=1,2$ are functions on $S^1$ valued in $ {\mathfrak h}$ which are over the base loop $q$ and satisfy $\pi^\# (\alpha_i)=u_i$. Also, the fact that $h_{2},u_{2}$ satisfies Condition (\ref{eq:tangentSpace}) has been exploited between the third and the fourth equality. In turn, this yields: 
\begin{eqnarray*}\omega (\delta_{a_1}, \delta_{a_2})  
& = &  - \int_{S^1} \left\langle k_2, v_1 \right\rangle  {\rm d}t-\int_{S^1}  \left\langle \pi^\# h_{1} , \alpha_{2}  \right\rangle {\rm d}t \\
 &+& \int_{S^1} \left\langle  \pi^\#\left(\frac{\diff \alpha_2}{\diff  t}\right)+  \frac{\diff \pi^\#}{\diff  t}(\alpha_2) - \frac{\partial{\pi^\#}}{{\partial q}} (\pi^\#(\alpha_2)) (p) , \alpha_{1}   \right\rangle   {\rm d}t  \\
& = & -\int_{S^1} \left\langle k_2, v_1 \right\rangle  {\rm d}t  -\int_{S^1}  \left\langle \pi^\# h_{1} , \alpha_{2}   \right\rangle {\rm d}t \\
&+& \int_{S^1} \left\langle  \pi^\#\left(\frac{\diff \alpha_2}{\diff  t}\right)+  \frac{\partial{\pi^\#}}{{\partial q}} (\pi^\# p) (\alpha_2) - \frac{\partial{\pi^\#}}{{\partial q}} (\pi^\#(\alpha_2)) (p) , \alpha_{1}   \right\rangle  {\rm d}t. 
 \end{eqnarray*}
 Re-arranging terms, by using the Jacobi identity in  the form of Lemma \ref{lemma:Jacobi} and the skew-symmetry of $\pi^\#$:
\begin{eqnarray*}
\omega (\delta_{a_1}, \delta_{a_2})& = &  -\int_{S^1} \langle k_2, v_1 \rangle  {\rm d}t -\int_{S^1}  \langle \pi^\# h_{1} , \alpha_{2}   \rangle {\rm d}t - \int_{S^1} \left\langle   \frac{\diff \alpha_2}{\diff  t}, \pi^\#\alpha_1 \right\rangle  {\rm d}t  \\
&+& \int_{S^1} \left\langle  \frac{\partial \pi^\#}{{\partial q}} (\pi^\# p) (\alpha_2) - \frac{\partial \pi^\#}{{\partial q}} (\pi^\#(\alpha_2)) (p) , 
\alpha_{1}   \right\rangle   {\rm d}t  \\
&=&  -\int_{S^1} \left\langle k_2, v_1 \right\rangle  {\rm d}t -\int_{S^1}  \left\langle \pi^\# h_{1} , \alpha_{2}   \right\rangle {\rm d}t
+ \int_{S^1} \left\langle   \alpha_2, \frac{\diff \pi^\#\alpha_1}{\diff  t} \right\rangle  {\rm d}t  \\
&+& \int_{S^1} \left\langle  \frac{\partial \pi^\#}{{\partial q}} (\pi^\#(\alpha_1)) (p), \alpha_{2} \right\rangle  {\rm d}t \\
 &=&  -\int_{S^1}  \left\langle  \pi^\# h_{1}  -  \frac{\diff u_1}{\diff  t} +  \frac{\partial \pi^\#}{{\partial q}} (u_1) (p), \alpha_{2}  \right\rangle {\rm d}t - \\
&-&  \int_{S^1} \langle k_2, v_1 \rangle  {\rm d}t. \\
  \end{eqnarray*}
This last identity proves both items in the result: first, the vanishing of $ \omega (\delta_{a_1}, \delta_{a_2})$ is automatic when 
 $h_{1},u_{1}$ satisfies Condition (\ref{eq:tangentSpace}) and $v_1=0$. Also, 
if $\omega (\delta_{a_1}, \delta_{a_2})$ vanishes for all $ \delta_{a_2} \in T_a C_\pi^{S^1}(S)$, it means that, since $k_2,\alpha_2$ can be chosen to be an arbitrary periodic  sections of $ {\mathfrak h}$
over the loop $q$ by Lemma \ref{lem:tgtconeIsSpace},  
it implies that  $h_{1},u_{1}$ satisfies Condition (\ref{eq:tangentSpace}) and $v_1=0$, i.e., in view of Lemma \ref{lem:tgtconeIsSpace} again, that $\delta_{a_1} \in T_a C_\pi^{S^1}(S)$. 
\end{proof}

\subsection{Periodic cotangent paths are not a submanifold}
\label{sec:counterExample}

In Theorem 1.5 in \cite{Cattaneo} and Theorem 1.3 in \cite{CattaneoFelderCoisotropic}, it is claimed that the result proved for $C_{\pi}(M) $, i.e. "it is a coisotropic submanifold of the Banach manifold of all $C^1$-paths valued in $T^{*}M$ if and only if $\pi$ is a Poisson bivector", still holds true for loops, i.e. for  $C_{\pi}^{S^1}(M)$. We are afraid that this generalization is problematic: it is certainly true that $C_{\pi}^{S^1}(M)$ is coisotropic in some sense, but it does not seem to true in general that it is Banach manifold.

\begin{counter-example}
There exists a Poisson manifold $(M,\pi)$ and a cotangent path on it such that the tangent cone $T_a C_{\pi}^{S^1}$  of cotangent loops is not a vector space. 
\end{counter-example}
This forbids cotangent loops from being a Banach submanifold (defined, in the Banach setting, as in Definition 3.2.1 of \cite{RJT} for instance - a definition large enough to include inverse image of 
regular values of smooth maps, which are themselves defined as being points where the differential of the submersion is onto and splits):
 \begin{lemma}
Let $C \subset{\mathfrak P} $ a Banach submanifold of a Banach manifold. 
For all $c \in C$, the tangent cone at $c$ is a vector space. 
\end{lemma}

We now construct our counter-example explicitly.
We will show that for $M = {\mathbb R}^{2}$, equipped with local coordinates $(x,y)$ and the following Poisson bivector:
\begin{equation}
\pi = x\pp{}{x}\wedge \pp{}{y},
\end{equation}
there exists a periodic cotangent path $c$ and two elements of the tangent cone $u,v \in T_c C_{\pi}^{S^1}(M)$ such that $u + v  \notin T_{c}C_{\pi}^{S^1}(M)$.

First, let us characterize periodic cotangent paths in this case. The equation defining cotangent paths, namely $q'= \pi^\#_{\gamma (t)}(p)$, reads:
\begin{equation}\label{equation1}
x' (t) = x (t) b (t) \hbox{ and } 
y' (t) = - x (t) a (t)
\end{equation}
with the understanding that $q(t)=(x(t),y(t)) $ and $p(t) = a(t) {\rm d}x + b(t) {\rm d} y $.

We consider the following cotangent path $c \in C_{\pi}^{S^1}(M)$, which is periodic because it is constant:
\begin{equation*}
  \begin{array}{cccc}
 c  :& S^{1}&\to& T^{*}{\mathbb R}^{2}\simeq {\mathbb R}^{2}\times {\mathbb R}^{2}\\
&t&\mapsto&(0,0)\times (0,0).
  \end{array}
  \end{equation*}%
We  now construct two elements in the tangent cone at $c$ of $ T_c C_{\pi}^{S^1}(M) $.
Consider the one-parameter family of periodic cotangent paths $\epsilon \to g_{\epsilon} $
and the one-parameter family of periodic cotangent paths  $\epsilon \to h_{\epsilon}$ given by:
$$ g_{\epsilon} = (\epsilon,0) \times (0,0) ~~\hbox{and}~~ h_{\epsilon} = (0,0) \times (0,\epsilon ).$$
For all fixed $\epsilon$, this is a cotangent path, which is periodic because it is constant. 
We set  $ u = \left. \frac{\diff g_(\epsilon )}{\diff \epsilon} \right|_{\epsilon=0}$ and $v = \left. \frac{\diff h_\epsilon }{\diff \epsilon} \right|_{\epsilon=0}$. 
Let us show that $u+v \not\in T_c C_{\pi}^{S^1}(M)$. Assume we are given a 1-parameter family of paths  $(x_{\epsilon }(t),y_{\epsilon }(t),a_{\epsilon }(t),b_{\epsilon }(t))$  valued in $C_{\pi}^{S^1}(M)$
whose derivative at $\epsilon=0$ is equal to
$u+v$, we have:
$$\left. \frac{\diff x_{\epsilon }(t)}{\diff \epsilon} \right|_{\epsilon=0}=1 ~\hbox{and}~ \left. \frac{\diff b_{\epsilon }(t)}{\diff \epsilon} \right|_{\epsilon=0}=1~ \hbox{and}~\left. \frac{\diff y_{\epsilon }(t)}{\diff \epsilon} \right|_{\epsilon=0}=0~\hbox{and}~\left. \frac{\diff a_{\epsilon }(t)}{\diff \epsilon} \right|_{\epsilon=0}=0$$
For all $\epsilon $ small enough and for all $t \in [0,1]$, the strict inequalities  $x_{\epsilon }(t)>0$ and $b_{\epsilon }(t)>0$ hold. 
Sense these $1$-parameter family of paths takes values in cotangent paths, by Equation (\ref{equation1}), we have:
$$\frac{\diff x_{\epsilon }(t)}{\diff t} =x_{\epsilon }(t)  b_{\epsilon }(t) >0.$$
But this forbids $x_\epsilon (t)$ from being a periodic function in $t$, since a periodic function can not be strictly increasing, hence the contradiction. 
We conclude that the path  $x_{\epsilon }(t)$ is not periodic, hence $u+v$ is not in the tangent cone at $c$.

\begin{remark}
 It is certainly true that cotangent loops defined over the regular part of a Poisson manifold form a Banach submanifold and that this submanifold is coisotropic.
 This implies that its closure (for the Banach topology) is a coisotropic set. But it deserves to be noticed that cotangent loops defined over the regular part 
 of a Poisson manifold are not a dense subset of cotangent loops, so it does not seem to be easy to prove (an equivalent in the Banach setting of) Theorem \ref{theo:2} with such a method.
\end{remark}

\end{document}